\documentclass[12pt]{amsart}
\usepackage{amsmath}
\usepackage{amssymb}
\usepackage{amscd}
\usepackage{graphicx}
\usepackage{calrsfs}
\usepackage{enumitem}
\usepackage{tikz-cd}
\usepackage{todonotes}
\usepackage[colorlinks=true,urlcolor=black, linkcolor=black, backref=page]{hyperref}   
\usepackage{etoolbox}
\usepackage{comment}
\usepackage{mathrsfs}

\newtheorem{thm}{Theorem}[section]
\newtheorem{cor}[thm]{Corollary}
\newtheorem{prop}[thm]{Proposition}

\newtheorem{quest}[thm]{Question}

\newtheorem*{openproblem*}{Problem}
\newtheorem*{quest*}{Question}
\newtheorem*{problem*}{Problem}

\theoremstyle{definition}

\theoremstyle{remark}
\newtheorem{rem}[thm]{Remark}

\newlength{\espaceavantspecialthm}
\newlength{\espaceapresspecialthm}
\newcommand{\R}{\mathbb{R}}
\newcommand{\N}{\mathbb{N}}

\newcommand{\Z}{\mathbb{Z}}

\def\diff{{\rm Diff}}

\def\S{{\bf S}}
\def\D{{\bf D}}
\addtocontents{toc}{\protect\setcounter{tocdepth}{1}}
\makeatletter
\let\c@equation\c@thm

\patchcmd{\subsection}{\bfseries}{\itshape}{}{}
\patchcmd{\@sect}{\@addpunct.}{}{}{}
\makeatother
\numberwithin{equation}{section}

\title[]{Problem List on Foliations and Diffeomorphism groups}

\begin{document}

\begin{abstract}
    This document collects problems proposed and discussed during the problem session at the conference  \href{https://conferences.cirm-math.fr/3082.html}{``Foliations and Diffeomorphism Groups"} (CIRM, 2024) organized by Hélène Eynard-Bontemps, Gaël Meigniez, Sam Nariman,
and Mehdi Yazdi. The problems were contributed by participants and have been lightly edited by the organizers for clarity and coherence.
\end{abstract}
\maketitle

\tableofcontents
\section{Introduction}
The study of diffeomorphism and homeomorphism groups of manifolds is deeply intertwined with foliation theory. Some of the most prominent open problems in the field, such as the Haefliger–Thurston conjecture and the L-space conjecture, highlight these connections.

Over the past decade, significant progress has been made on multiple fronts, including the dynamics of group actions on manifolds, the regularity of actions and foliations, algebraic properties of diffeomorphism groups, invariants of foliations, the Mather-Thurston homology equivalence, the group cohomology of diffeomorphism groups and their boundedness properties, the topology of spaces of foliations, representations of surface groups into diffeomorphism groups of the circle, and actions of 3-manifold groups on 1-manifolds via taut foliations.

The goal of this conference was to bring together a diverse group of experts and young researchers in foliation theory, diffeomorphism groups, 3-manifold topology, bounded cohomology, and 1-dimensional dynamics to exchange ideas and share insights. It is a long-standing tradition in conferences on foliation theory to compile a problem list, and we aimed to continue this practice. As organizers, we prepared a set of sample problems related to the conference themes and invited contributions from speakers and participants. During the conference, we slightly edited the submitted problems and incorporated additional ones that emerged from discussions in the problem session.

The problems vary in difficulty and background requirements. To keep the document concise, we did not add background material to shorter subsections but instead provided references and organized the problems thematically. Each subsection begins with the name of the person who contributed problems in that area.

We hope this problem list reflects the key questions that participants have been exploring and serves as an inspiration for others to engage with these topics.

\section{Taut foliations, pseudo-Anosov flows, and contact structures}

\subsection{Taut foliations and finite covers}  

\subsubsection*{Mehdi Yazdi}

Given a (transversely oriented) taut foliation $\mathcal{F}$ of a closed oriented 3-manifold, define the Euler class $e(\mathcal{F}) \in H^2(M ; \R)$ as the Euler class of the oriented tangent plane field to $\mathcal{F}$. Thurston \cite{thurston1986norm} showed that taut foliations of closed hyperbolic 3-manifolds have Euler classes of dual Thurston norm at most one. Conversely, he conjectured that any even integral class $a \in H^2(M ;\R)$ of dual norm one is the Euler class of a taut foliation on $M$. This was negatively answered by Gabai and Yazdi \cite{gabai2020fullymarked, yazdi2020thurston}. The following is a virtual version of Thurston's conjecture. 

\begin{quest}[Virtual Euler class one conjecture]
    Let $M$ be a closed hyperbolic 3-manifold whose first Betti number is at least one. Let $a \in H^2(M ; \R)$ be an (even) integral second cohomology class of dual Thurston norm one. Is there a finite cover $p \colon \Tilde{M} \rightarrow M$ and a taut foliation $\mathcal{F}$ of $\Tilde{M}$ such that the Euler class $e(\mathcal{F})$ of $\mathcal{F}$ is equal to $p^*(a)$? Here $p^* \colon H^1(M) \rightarrow H^1(\Tilde{M})$, and we have used Poincar\'{e} duality to identify $H^1(M) \cong H^2(M)$ and $H^1(\Tilde{M}) \cong H^2(\Tilde{M})$?
\end{quest}

If $a \in H^2(M ; \R)$ is a \emph{rational} class of norm one, Yi Liu \cite{liu2024criterion} has given a sufficient criterion in terms of the Alexander polynomial for $a$ to be virtually realised as the Euler class of a taut foliation on some finite cover $p \colon \Tilde{M} \rightarrow M$. Using his criterion, he has shown the following. 

\begin{thm}[\cite{liu2024criterion}]
For $b=2$ and for $b=3$, there exists some oriented closed hyperbolic
3–manifold M with $b_1(M) = b$, such that every rational point $a \in H^2(M; \R)$ of dual norm one is virtually realizable by a taut foliation. Indeed, the pullback of $a$ by some finite cyclic cover of $M$ dual to some primitive cohomology class in $H^1(M; \mathbb{Z})$ is the real Euler class of some transversely oriented taut foliation on the covering manifold.
\end{thm}

\subsection{Pseudo-Anosov flows transverse to taut foliations}

\subsubsection*{Jonathan Zung}
It is conjectured that any taut foliation on an atoroidal 3-manifold admits an (almost-)transverse pseudo-Anosov or periodic flow. It is natural to ask the converse question: given a pseudo-Anosov or periodic flow, when does there exist an (almost-) transverse foliation? The periodic/Seifert-fibered case is completely understood. The case where the base orbifold has one singular fiber is treated by the Milnor--Wood inequality. The case of small Seifert-fibered spaces is the subject of the Jankins--Neumann conjecture resolved by Naimi \cite{naimi1994foliations}, and also explained by Calegari--Walker \cite{calegari2011ziggurats}. The pseudo-Anosov case is wide open.

\begin{quest}
    Let $\varphi$ be a transitive pseudo-Anosov flow on a closed, oriented 3-manifold, and let $\gamma_1,\dots, \gamma_k$ be a collection of closed orbits. Let $Z\subseteq \R^k$ be the closure of the set of multislopes such that simultaneous Dehn surgery on the $\gamma_i$'s results in a pseudo-Anosov flow admitting a transverse codimension 1 foliation. What structure does the set $Z$ have?
\end{quest}

\begin{figure}
\includegraphics[width=2.5in]{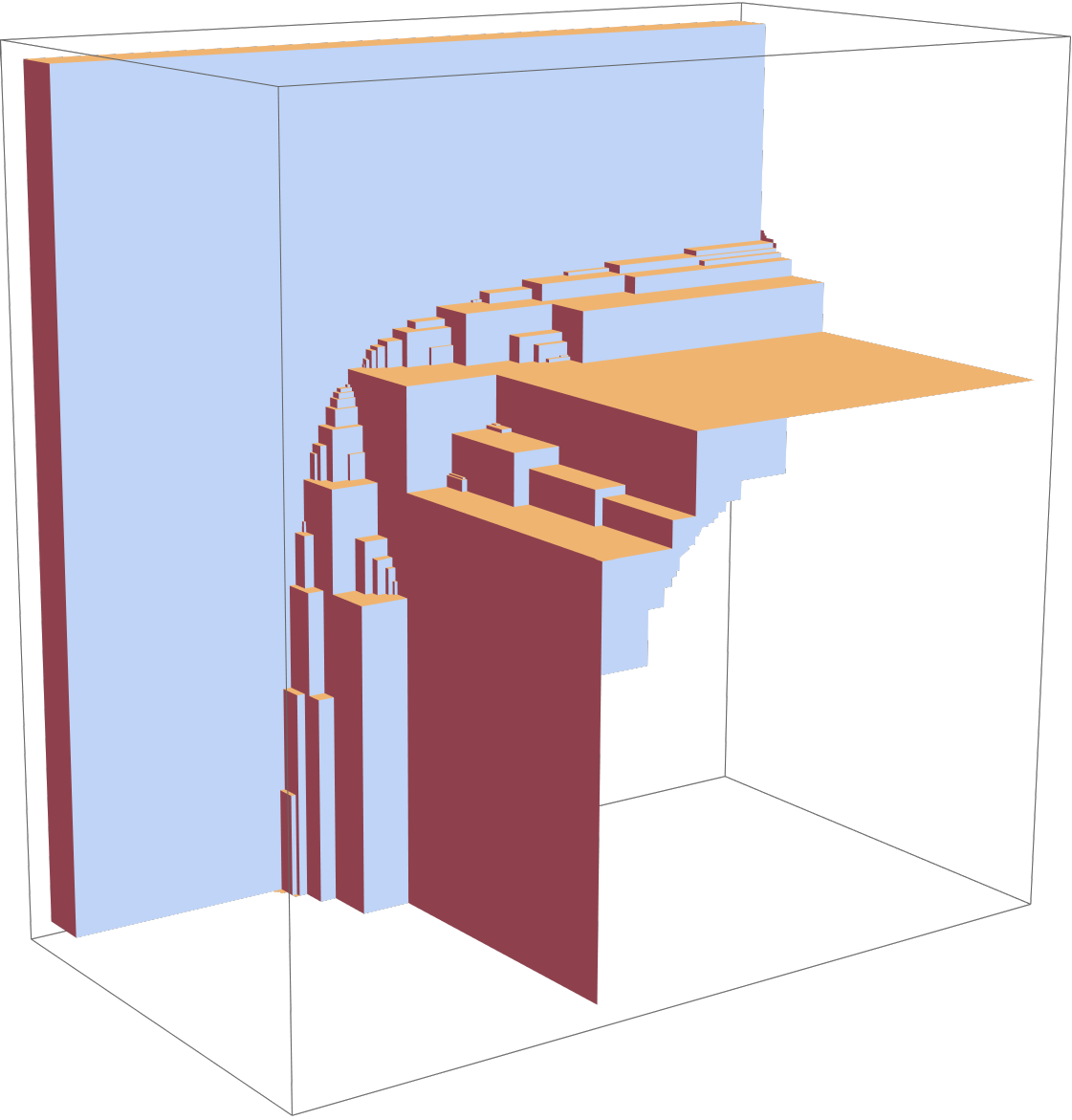}
\caption{An approximation to the ziggurat for Dehn surgeries on L6a5.}\label{fig:ziggurat}
\end{figure}

Here are some things we know about $Z$:
\begin{enumerate}
\item (Roberts) There are examples where $Z=\R^k$.
\item There are examples where $Z$ is a proper subset of $\R^k$. For example, some surgeries may be L-spaces which do not admit taut foliations.
\item (Massoni) Suppose $x,y,z\in \R^k$, and $x_i < y_i < z_i$ for all $1\leq i \leq k$. If $x,z\in Z$, then $y\in Z$.
\end{enumerate}
Experiments show that $Z$ has a stepwise structure; see Figure \ref{fig:ziggurat}. Following Calegari and Walker, who encountered similar sets in their study of maximal rotation numbers, we call such sets ``ziggurats''. Say that a set $S\subset \R^k$ is generated by a subset $T\subset S$ if $S$ is the smallest closed set containing $T$ and satisfying (3). We call the elements of a minimal generating set for $\overline Z$ the \emph{corners} of $\overline Z$. It appears that the closure of $Z$ is generated by a small number of corners. We expect that many of the properties that Calegari and Walker observed continue to hold in the pseudo-Anosov setting.

\begin{quest}
Establish the following properties of $Z$:
\begin{enumerate}
\item The corners of $Z$ are rational points.
\item The accumulation points of the corners are precisely the $\S^1\times \S^2$ surgeries.
\end{enumerate}

\end{quest}

\subsection{Contact pairs and taut foliations}

\subsubsection*{Thomas Massoni} In my talk, I discussed the following result: 
\begin{thm}[\cite{massoni2024taut}]
    Let $(\xi_-, \xi_+)$ be a pair of (co)oriented contact structures with opposite signs on a closed, connected, oriented, smooth $3$-manifold $M$. Assume that there exists a volume preserving vector field $Z$ which is positively transverse to both $\xi_-$ and $\xi_+$. Then $M$ admits a cooriented taut $C^0$-foliation.
\end{thm}

While this theorem already leads to some interesting applications (as discussed in my talk and in Jonathan Zung's talk), it would be convenient to relax the volume preserving condition on $Z$ and only impose ``intrinsic'' conditions on $\xi_-$ and $\xi_+$, i.e., conditions on the contact homotopy classes of $\xi_-$ and $\xi_+$, and not on the geometry of their intersection. The following question is a first step in this direction.

\begin{quest}
    Let $(\xi_-, \xi_+)$ be a pair of cooriented contact structures with opposite signs on $M$, and assume that there exists a vector field $Z$ transverse to both $\xi_-$ and $\xi_+$, and that both of $\xi_-$ and $\xi_+$ are tight. Does $M$ admit a cooriented \emph{Reebless} $C^0$-foliation?
\end{quest}

Refinements of the Eliashberg--Thurston theorem by Colin and Bowden imply that any Reebless aspherical $C^0$-foliation approximates to such pairs. As a warm-up, one could investigate the following question raised by Colin and Firmo in~\cite{colin2011paires}:

\begin{quest}
    Let $(\xi_-, \xi_+)$ be a tight positive contact pair as in the previous question. Is $M$ irreducible?
\end{quest}

A further step would be to remove the existence of $Z$ and only assume that $\xi_-$ and $\xi_+$ are homotopic as (cooriented) plane fields. To avoid complications, let us now assume that $M$ is an irreducible rational homology sphere.

\begin{quest}
    Let $(\xi_-, \xi_+)$ be a pair of (co)oriented contact structures with opposite signs on $M$ which are homotopic as plane fields, and assume that both of $\xi_-$ and $\xi_+$ are tight. Does $M$ admit a cooriented \emph{Reebless} $C^0$-foliation?
\end{quest}

We expect that there could be gauge theoretic obstructions preventing us to deform $(\xi_-, \xi_+)$ into a positive pair. It would be reasonable to add extra assumptions on the \emph{contact invariants} of $\xi_-$ and $\xi_+$ in Heegaard Floer homology. Recently, Francesco Lin~\cite{lin2023remark} showed that taut foliations approximate to contact structures whose contact invariants pair to $1$, for a natural duality pairing on Heegaard Floer homology. 

\begin{quest}
    Let $(\xi_-, \xi_+)$ be a pair of (co)oriented contact structures with opposite signs on $M$, and assume that their contact invariants pair to $1$: 
    \begin{align} \label{eq:pairing}
        \langle c(\xi_-), c(\xi_+) \rangle = 1.
    \end{align}
    Does $M$ admit a Reebless foliation?
\end{quest}

We note that the condition~\eqref{eq:pairing} is ``intrinsic'', and implies that $\xi_-$ and $\xi_+$ are tight and homotopic as plane fields. However, it is unclear how to interpret this condition geometrically.

\subsection{Euler classes of tight contact structures}

\subsubsection*{Mehdi Yazdi}

Eliashberg \cite{eliashberg1992contact} showed that if $\xi$ is a tight contact structure on a closed 3-manifold, then the Euler class $e(\xi)$ satisfies the inequality 
\[ \langle e(\xi), [S] \rangle \leq |\chi_-(S)|, \]
for all embedded surfaces $S$. In other words the dual Thurston norm of $\xi$ is at most one. This is the analogue of Thurston inequality \cite{thurston1986norm} for taut foliations. Every taut foliation on an irreducible 3-manifold can be $C^0$ approximated by tight contact structures \cite{eliashberg1998confoliations, bowden2016approximating, kazez2015approximating}, and so the set of Euler classes of tight contact structures on $M $ contains the set of Euler classes of taut foliations on $M$. This containement could be strict, for example in the case of the three-sphere. This motivates the following contact version of Thurston's conjecture. See \cite{sivek2023thurston}.
\begin{quest}
    Let $M$ be a closed hyperbolic 3-manifold whose first Betti number is at least one. Is every even integral second cohomology class $a \in H^2(M ; \R)$ of dual Thurston norm one realised as the Euler class of some tight contact structure on $M$?
\end{quest}

\subsection{Engel structures}
\subsubsection*{\'Alvaro del Pino}

An \textbf{Engel structure} is a maximally non-integrable plane field $\mathcal{D} \subset TM$ in a $4$-dimensional manifold $M$. More precisely:
\begin{itemize}
\item $\mathcal{D}$ is a rank-$2$ distribution,
\item it generates a hyperplane field $\mathcal{E} := [\mathcal{D},\mathcal{D}]$ by Lie bracket,
\item $\mathcal{E}$ is even-contact (i.e. satisfies $TM := [\mathcal{E},\mathcal{E}]$, meaning that it is maximally non-integrable itself).
\end{itemize}

We think of Engel structures as analogues of contact structures. Indeed, a generic plane field in dimension $4$, at a generic point, will be Engel. Moreover, Engel structures have a local model (they are locally isomorphic to the canonical distribution in the space of $2$-jets $J^2(\R,\R)$). The following result (attributed sometimes to Cartan) is explained in a nice paper of Richard Montgomery: the only distributions satisfying these two properties are the line fields, contact structures, even-contact structures, and Engel structures.

The group of Engel automorphisms of the local model is very rich: it is isomorphic to the group of contactomorphisms of standard tight $\R^3$. However, automorphisms do not globalize (we cannot ``cut them off''). In fact, it is known, due to work of R. Montgomery, that there are Engel manifolds with very small automorphism group. He posed the question of whether there exists an Engel manifold with trivial automorphism group. One such example was recently found by Koji Yamazaki, but the ambient manifold is non-compact. The following question remains open:
\begin{quest}
Is there a closed Engel manifold with trivial automorphism group?
\end{quest}
Yoshi Mitsumatsu has suggested a variation of this question:
\begin{quest}
Is the automorphism group of a generic Engel manifold trivial?
\end{quest}

Now to something different (but related). Engel manifolds have canonical dynamics. Namely, there exists a unique line field $\mathcal{W} \subset \mathcal{D}$, called the \textbf{kernel} or characteristic line field, with the property $[\mathcal{W},\mathcal{E}] = \mathcal{E}$. We think of it as follows: $\mathcal{E}$ is invariant under the holonomy of $\mathcal{W}$ (so it can be seen as a contact structure in its leaf space), while $\mathcal{D}$ rotates along $\mathcal{W}$ (i.e. we imagine it as rotating within the leaf space contact structure).

Nicola Pia has some work constraining the dynamics of the kernel line field $\mathcal{W}$, but overall we know very little. A crucial remark is that there is no stability for Engel structures, meaning that perturbing an Engel structure will produce some other Engel structure that in general won't be diffeomorphic to the original one. A vague question is then:
\begin{quest}
Are there any dynamical properties of $\mathcal{W}$ that are persistent under homotopies of the Engel structure?
\end{quest}

Here is a very concrete situation: Yoshi Mitsumatsu has a recipe, which he calls the \textbf{prequantum prolongation}, to produce Engel structures. The necessary input is a contact 3-manifold $(M,\xi)$, a line field $L$ tangent to $\xi$, and a non-vanishing closed 2-form $\omega$, of integral class, annihilating $L$. This implies that $L$ is spanned by a volume preserving vector field. The output Engel manifold is the $\mathbb{S}^1$-bundle over $M$ with curvature $\omega$, $\mathcal{E}$ is the corresponding connection, and $\mathcal{D}$ and $\mathcal{W}$ are the lifts of $\xi$ and $L$, respectively, to the connection. Mitsumatsu's paper contains very nice examples, but one could ask:
\begin{quest}
How does one systematically produce examples of a contact 3-manifold $(M,\xi)$, a line field $L$ tangent to $\xi$, and a non-vanishing closed 2-form $\omega$, of integral class, annihilating $L$?
\end{quest}
I did not think about this much, but perhaps one can drop the integrality condition by the following argument: if $\omega$ is not integral, we can perturb a bit and scale to make it so. This perturbs $L$. If this perturbation is small, it can probably be followed by a perturbation of $\xi$.

Another topic that is of interest is whether Engel structures interact with other geometric structures. Here are two questions:
\begin{quest}
Consider a $4$-manifold endowed with a foliation $\mathcal{F}$ by surfaces. Can the foliation be perturbed to an Engel structure?
\end{quest}
Do note that an Engel structure provides a splitting of TM into line bundles. It follows that we have to assume that $T\mathcal{F}$ and its complement split.

\begin{quest}
Consider a $4$-manifold endowed with a codimension-$1$ foliation $\mathcal{F}$. Suppose moreover that there is a plane field $\xi \subset T\mathcal{F}$ that is a leafwise contact structure. Assume that $\xi$ splits as a sum of line fields. Can this data be perturbed to an Engel structure?
\end{quest}
Casals, Presas and myself observed in a paper that the $h$-principle for overtwisted contact structures applies also in the foliated setting, so one can indeed construct many examples of such data.

\subsection{Foliation by planes of $\R^3$}

\subsubsection*{Takashi Tsuboi} I would like to conquer some doubt on Palmeira's theorem for foliation by planes $\mathcal{F}$ of $\R^3$.
In his paper there is a proposition saying that an infinite braid can be extended as an isotopy of plane, which has counter-examples. This is used to show that the leaf space determines the foliation. 
\begin{quest}
Find an alternative argument.
\end{quest}
I would like to explain a little more. In the paper \cite{palmeira1978open} by Carlos Frederico Borges Palmeira, he claims that for foliations $\mathcal{F}$ and $\mathcal{F}'$ by $(n-1)$-planes of simply-connected $n$-dimensional ($n\geqq 3$) manifolds $M$ and $M'$, respectively, if the quotient non-Hausdorff 1-dimensional manifold $M/\mathcal{F}$ and $M'/\mathcal{F}'$ are diffeomorphic then there is a foliation preserving diffeomorphism $(M,\mathcal{F})$ and $(M',\mathcal{F}')$ inducing the diffeomorphism in the quotient. In the proof of this claim he uses his main lemma on page 118, which is the case of the theorem for saturated closed subsets $N\subset M$ and $N' \subset M'$ such that the interiors of $N/\mathcal{F}$ and $N'/\mathcal{F}'$ are diffeomorphic to $\R$. In the argument in the proof of the main lemma, he takes a family of disjoint neighborhoods of leaves of the boundary and asserts that this ensures the construction of foliation preserving diffeomorphism. In a part of my master thesis of 1978, \emph{On foliations of $\R^n$ whose leaves are homeomorphic to $\R^{n-1}$} \cite{tsuboimastersthesis}, I proved a weaker result only for $n=4$ or $n\geqq 6$. I had a reason to exclude $n=3$. That is because a locally supported infinite braid of the plane may not be extended to an isotopy of $\R^2$. Such an example is illustrated in Figures 5 and 6 of my thesis. It might be possible to argue without thinking about such isotopies. In higher dimensions, this braid question does not appear by a general-position argument.

\section{Bounded cohomology}

\subsection{Bounded cohomology of diffeomorphism groups}
\subsubsection*{Mehdi Yazdi}
Let $\omega$ be a volume-form on $\S^2$ and $\mathrm{Diff}_0(\S^2, \mathrm{vol})$ be the group of volume-preserving diffeomorphisms of the 2-sphere that are in the same component as the identity map. The isomorphism type of the group $\mathrm{Diff}_0(\S^2, \mathrm{vol})$ does not depend on the choice of $\omega$. Banyaga \cite{banyaga1978structure} showed that $\mathrm{Diff}_0(\S^2, \mathrm{vol})$ is a simple group. Ghys and Gambaudo \cite{gambaudo2004commutators} showed that there are plenty of non-trivial quasimorphisms on the group $\mathrm{Diff}_0(\S^2, \mathrm{vol})$, in particular this group is not uniformly perfect. The following question is posed by Burago, Ivanov, and Poltervoich \cite{burago2008conjugation}.
\begin{quest}
    \begin{enumerate}
        \item Is $\mathrm{Diff}_0(\S^3, \mathrm{vol})$ uniformly perfect?
        \item Does $\mathrm{Diff}_0(\S^3, \mathrm{vol})$ admit any non-trivial quasimorphism?
    \end{enumerate}
    
\end{quest}

\subsubsection*{Sam Nariman} 
Let $\text{Diff}_0(\D^2, *)$ be the group of smooth orientation-preserving diffeomorphisms of $\D^2$ that fix a marked point $*$ say at the origin. There are two bounded Euler classes on this group; one is induced by the restriction to the boundary and the other is by taking derivative at the origin. The difference between these two classes gives a non-trivial quasimorphism on $\text{Diff}_0(\D^2, *)$.

\begin{quest}
Is the space of homogenous quasimorphisms on $\textnormal{Diff}_0(\D^2, *)$ one dimensional?
\end{quest}

If $A$ is a closed annulus, the space of homogeneous quasimorphisms on $\text{Diff}_0(A)$ is shown to be one dimensional by Militon \cite{militon2014commutator}. In fact the bounded cohomology of $\text{Diff}_0(A)$ is generated by two bounded Euler classes \cite{fournier2024bounded}, but the case of the disk with marked point is not known (at least to me).

\subsection{Bott Vanishing theorem in bounded cohomology}

\subsubsection*{Jonathan Bowden} Do we have the Bott vanishing theorem in bounded cohomology? More concretely.

\begin{quest}
    Using the Bott vanishing theorem, Morita showed that the MMM classes $\kappa_i$ for surface bundles vanish in $H^*(\textnormal{Diff}(\Sigma_g);\R)$ for $i>2$. Do they also vanish in bounded cohomology $H^*_b(\textnormal{Diff}(\Sigma_g); \R)$?
\end{quest}

\subsection{Simple orderable groups}

\subsubsection*{Xiaolei Wu}

By now, there are several constructions of simple orderable groups \cite{hyde2019finitely,hyde2023finitely, hyde2023two,bon2020groups}. However, only the examples in \cite{hyde2023finitely} of Hyde and Lodha are finitely presented.  In fact, they are even of type $F_\infty$. Recall a group $G$ is of type $F_n$ if there is a $K(G,1)$ with finite $n$-skeleton, and it is of type $F_\infty$ if it is of type $F_n$ for any $n$.

\begin{quest}
    For any $n\geq 2$, is there a simple orderable group of type $F_n$ but not $F_{n+1}$?
\end{quest}

If one drops the orderable condition, there are two constructions \cite{skipper2019simple, belk2022twisted}. A related question of Navas \cite[Question 8]{navas2018group} was the following: does there exist a finitely generated, left-orderable group $G$ without surjection onto $\mathbb{Z}$ and $H^2_b(G, \R) = 0$? The question was answered positively by Fournier-Facioa and  Lodha  in \cite{fournier2023second}. In fact, the simple orderable groups of Matte Bon and Triestino are even boundedly acyclic \cite[Corollary 1.16]{campagnolo2023algebraic}. But the only known finitely presented simple orderable groups are not boundedly acyclic  \cite[Corollary 1.7]{hyde2023finitely}. Thus we ask the following. 

\begin{quest}
    Is there a finitely presented, or even type $F_\infty$, boundedly acyclic simple orderable group?
\end{quest}

Another interesting question would be to find simple orderable groups of finite cohomological dimension. A candidate would be the Burger-Mozes groups \cite{burger2000lattices}.

\begin{quest}[Well-Known]
    Are Burger--Mozes' examples of simple groups orderable?
\end{quest}

%
%
%
%


\subsection{Milnor-Wood inequality, Euler class and Godbillon-Vey class}
\subsubsection*{Sam Nariman}

The boundedness of the Euler class for flat $\S^1$-bundles is known as the Milnor-Wood inequality \cite{milnor1958existence}, \cite{wood1971bundles}. For a $C^r$-flat $\S^1$-bundle $ E\xrightarrow{p} \Sigma_g$ over a closed oriented surface $\Sigma_g$ of genus $g>0$ states
\[
|\langle \mathscr{E}(p), [\Sigma_g]\rangle|\leq 2g-2,
\]
where $\mathscr{E}(p)$ is the Euler class of the bundle $E$. Freedman inspired by physics considered controlled Milnor-Wood inequality and controlled Mather-Thurston \cite{freedman2023controlled}. His controlled Milnor's inequality is as follows. He says ``Let  $g_{\mathscr{E}(p)}(\epsilon)$ be the smallest genus so that there is a generating set $S$ for $\pi_1(\Sigma_g)$ so that the circle bundle admits a topologically flat connection with representation $\rho\colon \pi_1(\Sigma_g) \to \text{PSL}(2,\R)$ so that elements in $\rho(S)$ are at most $\epsilon$ away in the sup norm from rotations $\text{U}(1)\subset \text{PSL}(2,\R)$. Define $g'_{\mathscr{E}(p)}(\epsilon)$ similarly by replacing $\text{PSL}(2,\R)$ with $\text{Homeo}_0(\S^1)$. Clearly $g'_{\mathscr{E}(p)}(\epsilon) \leq  g_{\mathscr{E}(p)}(\epsilon)$. It seems reasonable in light of \cite{wood1971bundles} to guess that they are actually equal."
\begin{thm}[Freedman]
The function $g_{\mathscr{E}(p)}(\epsilon)$, for $\epsilon>0$, is defined into the natural numbers union $0$. For $\epsilon$ sufficiently small it obeys the upper bound: 
\[
g_{\mathscr{E}(p)}(\epsilon) \leq  \frac{2\pi|\mathscr{E}(p)|}{\epsilon^2-O(\epsilon^3)}
\]
 where the error term satisfies $\frac{O(\epsilon^3)}{ \epsilon^3} < c$ , some constant  $c> 0$.
\end{thm}
\begin{quest}[Freedman]
Is that true that for small $\epsilon>0$, we have $g_{\mathscr{E}(p)}(\epsilon) \geq  \frac{c|\mathscr{E}(p)|}{\epsilon^2+O(\epsilon^3)}$, for some constant $c$ (perhaps $c=2\pi$)?
\end{quest}

Recall that the Euler class and Godbillon-Vey class are linearly dependent when the holonomy lands in $\text{PSL}(2,\R)$. But they are linearly independent classes in $H^2(\text{Diff}_0(\S^1);\R)$. We know that the GV class is unbounded but inspired by Freedman's theorem, we can ask the following.
\begin{quest}
Suppose $\epsilon>0$ is a small number. We have a representation  $\rho\colon \pi_1(\Sigma_g) \to \textnormal{Diff}_0(\S^1)$ such that for a generating set $S$ of $\pi_1(\Sigma_g)$, the set $\rho(S)$ are diffeomorphisms that are at most $\epsilon$ away from $\textnormal{PSL}_2(\R)$. Can we bound the GV class for the representation $\rho$ in terms of $g$  and $\epsilon$?
\end{quest}

To generalize the Milnor-Wood inequality to higher dimensions, Ghys \cite[S~F.1]{langevin1992list} posed a question of whether a generalization of Milnor--Wood inequality for flat $S^3$-bundles holds. Monod and I answered Ghys' question negatively by showing that 
\begin{thm}(\cite[Theorem 1.8]{monod2023bounded})
We have  $H^{4}_b(B\textnormal{Diff}^{r, \delta}_\circ(\S^3))=0$ for all $r\neq 4$. The Euler class and the first Pontryagin class $p_1$ in $H^4(B\textnormal{Diff}^{r, \delta}_\circ(\S^3))$ are unbounded.
\end{thm}
\begin{quest}
Is there a geometric proof of the above theorem? i.e. can one find a family of representations of a fundamental group of $4$-manifolds into $\textnormal{Homeo}_0(\S^3)$ that exhibits the unboundedness of the Euler class?
\end{quest}
Recently Zixiang Zhou posted a preprint in which he shows that $\text{Diff}(\S^n)$ is bounded acyclic for all $n\geq 4$! Therefore, the Euler class in $H^{2n}(B\text{Diff}^{\delta}(\S^{2n-1}))$ is not bounded as long as $n>1$. But the case of flat $\S^3$-bundles is particularly interesting because of its relation to Freedman's work on controlled Mather-Thurston theory. He asked the following question about Hopf bundle $\S^3\to \S^7\to \S^4$. 
\begin{quest}[Freedman]
Suppose $\S^3\to \S^7\to \S^4$ is a ``generalized'' Hopf fibration. Does there exist a homology $4$-sphere $H$ and a degree one map $H\to \S^4$ such that the pull-back of the generalized Hopf fibration is a flat $\S^3$-bundle?
\end{quest}
Recall that Mather-Thurston's theory for $C^0$-foliations implies that any $\S^3$-bundles $\S^3\to E\to M$ (in general manifold bundles) are bordant to a flat $C^0$ $\S^3$-bundle $\S^3\to E'\to M'$. Freedman's point of view is the bordism between $M$ and $M'$ that Mather-Thurston says it exists, should be semi-s-cobordism. In particular, we expect that the fundamental group of $M'$ should be more complicated. One naive guess of the meaning of the unboundedness of the Euler class and $p_1$ could be there should be a version of Mather-Thurston's theorem for $C^0$-foliations that bound the simplicial volume of the base. Roughly speaking, in the case of $\S^3$-bundles, the following might be reasonable to expect.
\begin{quest}
Suppose we are given an $\S^3$-bundle $\S^3\to E\xrightarrow{p} M$ over a $4$-manifold $M$ with nontrivial Euler class and $p_1$. Mather-Thurston's theory says we can find a bordism of this bundle to a $C^0$-flat bundle $\S^3\to E'\to M'$. Can one bound the simplicial volume of $M'$ so that it does not depend on the bundle map $p$?
\end{quest}

\subsection{Milnor-Wood inequality for 3-manifolds}

\subsubsection*{Mehdi Yazdi}

If $M$ is a manifold and $\rho \colon \pi_1(M) \rightarrow \mathrm{Homeo}_+(\S^1)$ is a homomorphism, then we can define $e(\rho) \in H^2(M ; \mathbb{Z})$ as the Euler class of the associated $\S^1$-bundle over $M$. By the Milnor-Wood inequality, the Euler class $e(\rho)$ is bounded in the sense that for every embedded surface $S$ in $M$
\[ \langle e(\rho) , [S] \rangle \leq |\chi(S)|. \]
If $M$ is a closed orientable 3-manifold, then the Milnor-Wood inequality can be restated as: the dual Thurston norm of $e(\rho)$ is at most one. The following is inspired by Thurston's Euler class one conjecture \cite{thurston1986norm, yazdi2020thurston} for taut foliations, and the universal circle construction of Thurston \cite{thurston1972foliations}, and Calegari and Dunfield \cite{calegari2003laminations}. In particular, if $\mathcal{F}$ is a taut foliation of $M$ and $e(\mathcal{F}) \in H^2(M ; \mathbb{Z})$ is the Euler class of the oriented tangent plane field to $\mathcal{F}$, then there is a faithful action $\rho \colon \pi_1(M) \rightarrow \mathrm{Homeo}_+(\S^1)$ (called the universal circle action) such that $e(\rho) = e(\mathcal{F}) \in H^2(M ; \mathbb{Z})$.

\begin{quest}
    Let $M$ be a closed hyperbolic 3-manifold whose first Betti number is at least one. Is every integral second cohomology class $a \in H^2(M ; \R)$ of dual Thurston norm one realised as the Euler class of some representation $\rho \colon \pi_1(M) \rightarrow \mathrm{Homeo}_+(\S^1)$?
\end{quest}

We can ask for other versions, where we consider cohomology classes of dual Thurston norm at most one, or we ask for faithful representations.

\subsection{Amenable actions and planar groups and the Euler class}
\subsubsection*{Sam Nariman} Hirsch and Thurston \cite{hirsch1975foliated} proved that the Euler class of amenable actions on spheres $\S^{2n-1}$ should vanish. This is interesting because now using bounded cohomology techniques, we know a non-constructive proof that the Euler class for flat $\S^{2n-1}$ bundles when $n>1$ is unbounded. However, by Hirsch-Thurston, it vanishes when the holonomy is amenable. 
\begin{quest}
How can one explain this phenomenon that the Euler class in $H^{2n}(\textnormal{Homeo}_0(\S^{2n-1}))$ is unbounded but it vanishes on amenable subgroups? Is it a ``hyperbolic" class in the sense of Gromov (see \cite{brunnbauer2024atoroidal})? Could it be that there is non-trivial bounded Euler class in $H^{2n}(\textnormal{Homeo}^{\text{vol}}_0(\S^{2n-1}); \R)$ where $\textnormal{Homeo}^{\text{vol}}_0(\S^{2n-1})$ means the volume preserving homeomorphisms of the sphere $\S^{2n-1}$?
\end{quest}

For Hirsch--Thurston's theorem, it is important for the fiber to be compact (in our case, they were spheres). Less is known when the fiber is not compact. Calegari  \cite{calegari2004circular} showed that the Euler class is unbounded in $H^2(\text{Homeo}_0(\R^2))$. So, surface group actions on $\R^2$ do not satisfy the Milnor-Wood inequality. However, he showed that for $\mathbb{Z}^2$ actions on $\R^2$ as long as the action is $C^1$, the Euler class has to vanish.

In a recent work with Fournier and Monod \cite{fournier2024bounded}, we showed that in fact $\text{Homeo}_0(\R^n)$ is bounded acyclic for all $n$. Therefore, in particular, the Euler class for actions on $\R^{2n}$ is unbounded. But we do not know if the analog of Hirsch--Thurston still holds in these cases.
\begin{quest}
Suppose $G$ is an amenable group acting on $\R^{2n}$ such that the action is at least $C^1$. Does the Euler class of such action vanish?
\end{quest}

Calegari's paper has many interesting points about the Euler class and planar group actions, but here I want to emphasize one of his questions:
\begin{quest}[Calegari]
Do analytic actions of surface groups on $\R^2$ satisfy the Milnor--Wood inequality?
\end{quest}
to which I also want to add
\begin{quest}
Do volume preserving actions of surface groups on $\R^2$ satisfy the Milnor--Wood inequality?
\end{quest}

The study of planar groups and amenable group actions on $\R^2$ could also lead to the generalization of Thurston's conjecture on $C^0$-stability of foliations with amenable holonomy. Recall that Thurston \cite{thurston1974generalization} generalized Reeb stability for $C^1$-foliations with compact leaves. He gave a counterexample of his generalization for $C^0$-codimension one foliations and conjectured that if the fundamental group of a compact leaf $L$ in a codimension-one transversely orientable foliation is amenable and if the first cohomology group $H^1(L; \R)$ is trivial, then $L$ has a neighborhood foliated as a product. This was later proved as a consequence of Witte-Morris’ theorem \cite{morris2006amenable} on the local indicability of amenable left orderable groups and Navas’ theorem \cite[Prop. 3]{mann2015left} on the left orderability of the group of germs of orientation-preserving homeomorphisms of the real line at the origin.

Thurston's stability hold for $C^1$-transversally oriented foliations in all codimensions. It also holds for $C^0$ codimension one foliations if we further assume that the fundamental group of the leaf is amenable, by the proof of Thurston's conjecture. 
\begin{quest}
Under what condition on the holonomy of a compact leaf $L$ in a $C^0$ transversally oriented foliation of codimension $>1$, the leaf $L$ is stable? For example, suppose $\mathcal{F}$ is a $C^0$ codimension $2$ foliation and $L$ is a compact leaf whose normal microbundle has a trivial Euler class, its fundamental group is amenable and $H^1(L; \R)=0$. Is $L$ stable?
\end{quest}

\subsection{Rigidity of group actions on the circle}

\subsubsection*{Hiraku Nozawa}

This is a question on a generalization of Mann's rigidity theorem for surface group actions on the circle to the case where the surface has cusps.

Let $\Sigma$ be a closed Riemann surface of genus $g>1$ and $\Gamma = \pi_1\Sigma$. For a $\Gamma$-action on $\S^1$, the Euler number $e(\rho) \in \mathbb{Z}$ is defined as the Euler number of the suspension bundle. The Milnor--Wood inequality 
\[
|e(\rho)| \leq |e(\Sigma)|
\]
says that the Euler number cannot be greater than that of the Fuchsian action
\[
\rho_0 : \Gamma \to \operatorname{Isom}_+ \D^2 \to \operatorname{Homeo}_+\partial \D^2,
\]
which comes from the uniformization of $\Sigma$.
 Matsumoto's rigidity theorem \cite{matsumoto1987some} says that the actions that attain the equality in the Milnor--Wood inequality are only the Fuchsian action up to semi-conjugacy.

\begin{thm}[Matsumoto] Let $\Sigma$ be a closed Riemann surface of genus $>1$ and $\Gamma=\pi_1\Sigma$. Let $\rho$ be a $\Gamma$-action on $\S^1$. If $e(\rho) = e(\Sigma)$, then $\rho$ is semi-conjugate to the Fuchsian action $\rho_0$, 
    i.e.\ $\exists \psi : \S^1 \to \S^1$ s.t.\  continuous of degree one, monotone and $(\rho, \rho_0)$-equivariant.
    i.e.,
    \[
\psi \circ \rho(\gamma) = \rho_0(\gamma) \circ \psi \quad (\forall \gamma \in \Gamma).
\]
\end{thm}

Mann \cite{mann2015spaces} generalized Matsumoto's rigidity to the \emph{geometric} actions: A $\Gamma$-action on $\S^1$ is called geometric if it is faithful and has a discrete image contained in a finite dimensional Lie group $G$ in $\operatorname{Homeo}_+\S^1$.

\begin{thm}[Mann]\label{thm:Mann} Let $\Sigma$ be a closed Riemann surface of genus $>1$ and $\Gamma=\pi_1\Sigma$. Let $\rho$ be a geometric $\Gamma$-action on $\S^1$. Let $C_\rho$ be the connected component of $\operatorname{Hom}(\Gamma,\operatorname{Homeo}_+\S^1)$ that contains $\rho$. Then every $\rho'$ in $C_\rho$ is semi-conjugate to $\rho$.
\end{thm}

Now let us assume that the Riemann surface $\Sigma$ may have a finite number of cusps. Burger--Iozzi--Wienhard generalized the Milnor--Wood inequality and Matsumoto's rigidity to Riemann surfaces with cusps. We should take care that the usual Euler number becomes trivial for $\Gamma$-actions, because the suspension bundle over $\Sigma$ is topologically trivial. Burger--Iozzi--Wienhard \cite{burger2010surface} defined the Euler number by using real coefficient bounded cohomology based on the work of Ghys. This Euler number is not homotopy invariant and is not an integer in general, but it is expressed by using Poincar\'{e}'s translation number as in the case of the classical Euler number.

The question in this section is the following:
\begin{quest}
Is it possible to generalize Mann's theorem (Theorem \ref{thm:Mann}) to the case where the surface has a finite number of cusps?    
\end{quest}

We should note that a part of the formulation of the theorem should be appropriately changed to be generalized. For example, 
the fundamental group of a surface with cusps is a free group. Thus, the entire space of actions on $\S^1$ is connected.

\section{Stability of foliations} 

\subsection{Using commuting vector fields}
\subsubsection*{Christian Bonatti}
Given an isolated fixed point of a diffeomorphism or an isolated zero of a vector field, the Poincar\'e--Hopf index allows us to decide if this point will persist under small perturbation of the diffeomorphism or of the vector field. On a compact manifold, the sum of the indices of singular points of a vector field is the Euler characteristic of the manifold.  This also holds for diffeomorphisms isotopic to the identity, and in general the sum is given by Poincar\' e Lefschetz formula.  The question I am formulating below comes from my interrogations about what could be an equivalent to that index for compact leaves of foliations.  For compact leaves of $1$-dimensional oriented foliations, the solution has been given by Fuller \cite{fuller1967index}. The fiber bundle with fiber the circle can be seen as the equivalent of the identity map or the zero vector field, and the Fuller index decides if foliations close to oriented circle bundles still have compact leaves close to a fiber. In fact, in that oriented setting, Fuller index is just the Euler characteristic of the basis. In the non-oriented case this sum is the class modulo $2$ of the Euler characteristic. 

In order to extend these results to compact leaves with an abelian fundamental group, I tried for a long time to answer the following question
\begin{quest} Consider a foliation $\mathcal{F}$ defined by a fibration with closed basis $B$ and fiber $F$ and $\pi_1(F)=\mathbb{Z}^k$.  Under what hypothesis does any foliation $\mathcal{G}$ close to $\mathcal{F}$ have a compact leaf close to a fiber? 
\end{quest}
This question is completely understood when $B=\S^1$ and the answer depends on the monodromy of the fibration, in particular on the existence of real positive eigenvalues of this monodromy. See in particular my paper with Sebasti\~ao Firmo, where we define a notion of index for codimension $1$ compact leaves with abelian fundamental group, and whose sum is an invariant depending on the monodromy and vanishes if and only if one may destroy every leaf of the fibrations \cite{bonatti1994feuilles}. 

The next step should be to consider torus-fibration with $2$-dimensional basis, and to decide, in terms of the monodromy of the fibration, if all fibers may be destroyed by small perturbation of the foliation, or if some of them need to persist. In my papers \cite{bonatti1989point} and \cite{bonatti1990stabilite}, I showed that the answer is positive when the monodromy is trivial and the basis has non trivial Euler characteristic, but I also built a specific monodromy for which the answer was negative. 

Thus a natural question, for understanding what could be the index I am looking for is
\begin{quest} Let $\mathcal{F}$ be the foliation defined by a fibration $\pi\colon M\to S$ where $S$ is a genus $2$ closed oriented surfaces and the fiber is the torus $T^2$.  Under what hypothesis on the monodromy does any foliation $\mathcal{G}$ close to $\mathcal{F}$ have a compact leaf close to a fiber? 
\end{quest}

The holonomy of any foliation $\mathcal{G}$ $C^1$-close to $\mathcal{F}$ consists locally of two commuting diffeomorphisms of the base surface $S$, for a given basis of the fundamental group of the fiber. However, due to the monodromy, one cannot choose globally such a basis. Therefore the holonomy of $\mathcal{G}$ looks like commuting diffeomorphisms statisfying some equivariance relation related with the monodromy.  A simpler way to build such commuting diffeomorphism is to build commuting vector fields, with a corresponding equivariance relation. This motivates the following problem:

Let S be a closed orientable surface of genus greater than $1$.
Let $ \rho : \pi_1(S) \to \mathrm{SL}(2,\mathbb{Z})$  be a representation.
 Let $\tilde S$ denote the universal cover of S and 
$\varphi :\pi_1(S)\to \text{Diff}(\tilde S)$ be the natural action by deck transformations.

One says that two commuting vector fields $(U,V)$ on $\tilde S$ are $\rho$-equivariant if
for any $p$  and $\gamma\in\pi_1(S)$ one has
\[ (U(\varphi(\gamma) (p)) , V(\varphi(\gamma)(p))= \rho(\gamma)(U(p),V(p)), \]
where  $A(U,V) =(aU+bV, cU+dV))$  if $A$ is the matrix 
$$\left(
\begin{array}{cc}
a & b\\
c & d
\end{array}
\right).$$
\begin{quest}
Characterize the representations $\rho$ for which any $\rho$-equivariant commuting vector fields $(U,V)$ have a common zero.
\end{quest}

Elon Lima proved that commuting vector fields for the trivial representation have a common zero.
I build an example of $C^1$ commuting pair of vector fields with no common zero, equivariant by $\rho$ if  $\rho(\pi_1(S))$ is generated by a unique hyperbolic matrix.

This question is essentially equivalent to the following problem:
\begin{quest}
Consider a torus $T^2$ bundle over S.  Under what hypothesis on the monodromy $\rho$ can we find a foliation with no compact leaves close to a fiber?
\end{quest}

My results \cite{bonatti1990stabilite} provide an example of two torus fibrations over the same closed surface $S$: the trivial and the one with a cyclic hyperbolic monodromy.  Let $F$ and $G$ be the foliations defined by these fibrations (the leaves are the fibers): for the trivial fibration any perturbation of the foliation $F$ \emph{preserves a fiber}. In contrast, there are perturbations opening every fiber of $G$. On the other hand $F$ and $G$ admit perturbations $\tilde F,\tilde G$ having finitely many compact leaves, and $\tilde F$ and $\tilde G$ are conjugated in a neighborhood of their compact leaves: thus, if there exists an index theory for compact leaves (as the Poincar\'e--Hopf index for vector fields) they should have the same index, if this index \emph{is a number}. But the sum of the indices for $F$ should  $\neq 0$ when the one for $G$ should vanish. My guess is that such an index should live in some homology class related to the fibrations, and the fact that one fibration is trivial and the other is twisted could solve the paradox. A solution of the questions could lead to an intuition of where the \emph{index} should live.

\subsection{Using  perturbations of group actions on the interval}

\subsubsection*{Christian Bonatti}
\begin{quest}Prove  every leaf of a $C^r$-generic codimension $1$ oriented and co-oriented foliation on a $3$-manifold  has genus $<2$.
\end{quest}
The difficulty of this question depends strongly on the regularity $r$ of the foliations and of the topology on the space of foliations. The case of compact leaves has been done in my paper with Firmo \cite{bonatti1994feuilles} for $r=1$ and by Tsuboi \cite{tsuboi1994hyperbolic} for $r\geq 2$. 

For noncompact leaves, the result could be implied by a positive answer to the following question.

\begin{quest} Consider four local $C^r$ diffeomorphisms $f_a,f_b,f_c,f_d$ of $R$ whose definition domain contains  $[-1,1]$ and having $0$ as a common fixed point. Let $h=[f_a,f_b][f_c,f_d]$ be the product of the commutators. 

Does there exist $\tilde f_a,\tilde f_b, \tilde f_c,\tilde f_d$ coinciding with $f_a,f_b,f_c,f_d$ out of $[-1,1]$, $C^r$-close to $f_a,f_b,f_c,f_d$, so that 
$$h=[\tilde f_a\tilde f_b][\tilde f_c\tilde f_d]$$
and so that $$Fix f_a\cap Fix f_b\cap Fix f_c \cap Fix f_d \cap (-1,1)=\emptyset.$$ 
\end{quest}
The case $h=id$ is the one we solved and which gives the answer for compact leaves. 

\section{Commuting diffeomorphisms}
\subsection{Topology of the space of commuting homeomorphisms}

\subsubsection*{Sam Nariman}

We learned from Andres' talks that the space of commuting elements in $\text{Homeo}_0(\S^1)$ is connected. The inclusion $$\text{Homeo}_0(\S^1)\hookrightarrow \text{Homeo}_0(\R^2)$$ is a weak homotopy equivalence. However, the space $\text{Hom}(\mathbb{Z}^2, \text{Homeo}_0(\R^2))$  is not connected since $\mathbb{Z}^2$-actions on $\R^2$ can have different Euler numbers.
\begin{quest}
Let $\textnormal{Hom}_k(\mathbb{Z}^2, \textnormal{Homeo}_0(\R^2))$ be the subspace of $\mathbb{Z}^2$-actions on $\R^2$   whose Euler number is $k$. Are subspaces $\textnormal{Hom}_k(\mathbb{Z}^2, \textnormal{Homeo}_0(\R^2))$ connected for each integer $k$?
\end{quest}

\subsection{2-cycles defined by commuting diffeomorphisms}
\subsubsection*{Takashi Tsuboi}
Since our knowledge of the commuting diffeomorphisms of the circle has deepened by the works of Navas, Eynard-Bontemps and others, I would like to pose the following question again 40 years after my paper.
\begin{quest}\label{t}
For commuting diffeomorphisms $f$ and $g$ of $\S^1$, to what extent we can say that 
the 2 cycle $(f,g)-(g,f)$ of $B\mathrm{Diff}(\S^1)^\delta$ is a boundary ?
\end{quest}
If $f$ and $g$ generate an infinite cyclic group, then it is trivially a boundary.
If  $f$ and $g$  belong to a flow on $\S^1$, then one can show that it is a boundary \cite{tsuboi1981cycles}.
\newline
\begin{rem}[Sam Nariman]
    One of the long-standing questions about codimension one foliations is whether the vanishing of the Godbillon-Vey class for a codimension one foliation $\mathcal{F}$ implies that $\mathcal{F}$ is foliated null-cobordant. The dynamical consequences of the vanishing of the Godbillon-Vey class have been extensively studied (see \cite{ghys1989invariant} and the introduction of \cite{tsuboi1981cycles}). In particular, it was proved by Herman that the Godbillon-Vey class of all flat $\S^1$-bundles over the $2$-torus $T^2$ vanishes. But it is still not known whether such foliations are foliated nullbordant, which motivates Question \ref{t}.
\end{rem} 

\subsection{The critical regularity of the centralizer}
\subsubsection*{Sang-hyun Kim}
Let $f$ be a $C^1$--diffeomorphism of a closed connected smooth manifold $M$. One may define its critical regularity
\[
\operatorname{CritReg}(f):=\sup\{ r\in[1,\infty)\mid f\text{ is }C^r\}\in[1,\infty].\]
It is standard to construct a diffeomorphism $f$ with 
$\operatorname{CritReg}(f)=r$ for each real $r\in[1,\infty]$; see \cite{kim2020diffeomorphism, mann2019reconstructing} for $\dim M=1$ and see~\cite{harrison1975unsmoothable,harrison1979unsmoothable} for $\dim M>1$. In fact, one may even require that no topological conjugacy of $f$ can increase the regularity when $\dim M>1$. Similarly, one can define
\begin{align*}
&\operatorname{CritCent}(f):=
\\
&\sup\{ r\in[1,\infty)\mid \text{ for some nontrivial }g\in\operatorname{Diff}^r(M)\text{ centralizes }f\}\\
&\ge\operatorname{CritReg}(f).
\end{align*}
Generically in the $C^1$--topology, the centralizer group of $f$ coincides with the cyclic group generated by $f$ (see~\cite{bonatti2009c}) and in this case, we have $\operatorname{CritCent}(f)=\operatorname{CritReg}(f)$.

\begin{quest}
For which $f\in\operatorname{Diff}^1(M)$, does one have \[\operatorname{CritCent}(f)>\operatorname{CritReg}(f)?\]
What are the possible values of $\operatorname{CritCent}(f)$ in this case?
\end{quest}

\section{Actions on 1-manifolds preserving a foliated structure}

\subsection{Actions on  the circle and bi-foliated planes}
\subsubsection*{Christian Bonatti}

Consider a group action $\Phi: G\to \text{Homeo}^+(\R^2)$ on the plane preserving 2 transverse (maybe with prong singularities) foliations $\mathcal{F}^+$ and $\mathcal{F}^-$.
This action always extends \cite{bonatti2024action} to \emph{the circle at infinity} associated to the pair of foliations $(\mathcal{F}^+,\mathcal{F}^-)$, as an action $\varphi: G\to Homeo^+(\S^1)$ called \emph{the action induced by $\Phi$}.
\begin{quest}
Can we see on the induced action $\varphi$ if $\Phi$ is smoothable or not?
\end{quest}

We can reformulate the question as follows:  
\begin{quest}Is there an action $\varphi$ on $\S^1$ which is induced by a topological action $\Phi$, but is not induced by any smooth action $\Psi$ on $\R^2$ preserving two transverse foliations? 
\end{quest}

The answer is \emph{no} for actions coming from transitive Anosov flows due to Shannon result \cite{shannon2020dehn} in his PhD thesis (topological transitive Anosov flows are equivalent to smooth Anosov flows), but it remains open for nontransitive topological Anosov flows. 

\subsubsection*{Kathryn Mann}
Here is a related question about bifoliated planes.  As explained in the article of Haefliger and Reeb \cite{haefliger1957varietes}, any (possibly non-Hausdorff) simply connected 1-manifold is the leaf space of a foliation of the plane.  
\begin{quest}
    Which pairs $\Lambda, \Lambda'$ of (possibly non-Hausdorff) simply connected 1-manifolds can be realized as the leaf spaces of a pair of transverse foliations of the plane?   What if you replace "pair" by    ``triple"?   

    One can also specialize this to the pairs of stable/unstable leaf spaces for Anosov flows, it is as yet not even known whether the stable leaf space determines the unstable.  
\end{quest}

\subsection{Amenable orderable groups}
\subsubsection*{Michele Triestino}

A well-known result of Witte Morris states that any amenable (left-)orderable group is locally indicable: any finitely generated non-trivial subgroup surjects to $\mathbb Z$. This characterizes amenable orderable groups among amenable groups, since locally indicable groups are always orderable.

Assume that $G$ is a finitely generated amenable orderable group. Given a left-order on $G$, the so-called \emph{dynamical realization} gives a faithful action of $G$ on $\R$ by orientation-preserving homeomorphisms. We propose to investigate the possible restrictions for the dynamics of such an action. More precisely, we propose to describe the possible actions of $G$ on $\R$, without global fixed points, up to \emph{semi-conjugacy}. Recall that two actions $\varphi_1,\varphi_2\colon G\to \mathrm{Homeo}_+\R$ without global fixed points are semi-conjugate if there exists a $(\varphi_1,\varphi_2)$-equivariant non-decreasing map $h\colon \R\to \R$.

The result of Witte Morris implies that $G$ always admits a non-trivial action by translations on $\R$. The possible actions by translations of $G$, up to (semi-)conjugacy, are in correspondence with the sphere $H^1(G;\R)\setminus \{0\}/(\lambda>0)$ of non-trivial morphisms $G\to (\R,+)$ up to positive rescaling.

When $G$ is of sub-exponential growth, this is a complete description, by a result of Plante \cite{plante1975foliations}. When $G$ is virtually solvable, it has been proved in \cite[Theorem 8.3.8]{brum2021locally} that any action without global fixed points is either by affine transformations, or \emph{laminar}: that is, one can find an interval $I\subset \R$ such that any two images are either nested or disjoint, and the union of all images covers $\R$.

The precise question we ask here is to exhibit examples of actions without global fixed points which are not of this nature:

\begin{quest}
    Let $G$ be a finitely generated amenable orderable group, and $\varphi\colon G\to \mathrm{Homeo}_+\R$ an action of $G$ without global fixed points. Is it true that $\varphi$ is either semi-conjugate to an affine action or laminar?
\end{quest}

The way we formulate this question is extremely optimistic, but in fact we have no idea of whether this should be true: as a big \textit{caveat}, notice that the standard action of Thompson's group $F$ on $(0,1)\cong \R$ is not of this nature.

\subsection{Laminar groups}

\subsubsection*{Harry Hyungryul Baik}
Above Michele defined a laminar group action on $\R$. One can also define a laminar group action on $\S^1$ by first defining a lamination on $\S^1$ as a closed subset of 
\[ ( \S^1 \times \S^1 - \Delta ) / (x, y) \sim (y, x)\] 
whose elements are pairwise unlinked. Here $\Delta$ is the diagonal \linebreak $\{ (x,x) | x \in \S^1\}$. Barthelm\'e-Bonatti-Mann defined a prelamination in the same way, except without imposing the closedness condition. (almost all) Hyperbolic surface groups are laminar as there are many geodesic laminations on any surface other than the three-punctured sphere and many hyperbolic 3-manifold groups are also laminar following work of Thurston, Calegari-Dunfield, Calegari, Fenley, and Frankel-Schleimer-Segerman. Hyperbolic surface groups are completely characterized as laminar groups in terms of their invariant laminations by myself, which later generalized to hyperbolic 2-orbifolds by myself and KyeongRo Kim. More recently, together with Hongtaek Jung and KyeongRo Kim, I showed that if a group acts on the circle with so-called a veering pair of invariant laminations, then the group is a 3-orbifold group \cite{baik2025groups}. If the action is cofinite (equivalently the resulting orbifold is compact) then the resulting orbifold is hyperbolic due to geometrization. Mann conjectured that if a group admits an Anosov-like action on a bifoliated plane, then the group is the fundamental group of a 3-manifold with a pseudo-Anosov flow and the bifoliated plane is the orbit space of the flow (we will refer it as the Anosov-like action conjecture). In the setting of Baik-Jung-Kim , if the action is cofinite and the veering pair of laminations has no crown gaps, the action induces an Anosov-like action. In this case, the 3-manifold can be constructed, and by the functoriality of the work of Frankel-Schleimer-Segerman, the induced bifoliated plane is the orbit space of a pseudo-Anosov flow without perfect fits, and also all such actions arise in this form. Hence, Baik-Jung-Kim \cite{baik2025groups} settles the Anosov-like action conjecture in the case of pseudo-Anosov flows without perfect fits, although the general case should be much more complicated. Via the work of Barthelm\'e-Bonatti-Mann, the viewpoint of the work of Baik-Jung-Kim can translate into the following question (I will write it in a bit informal way). 

\begin{quest} 
Which pair of prelaminations forces any group action that preserves it to induce an Ansov-like action on the corresponding bifoliated plane? 
\end{quest}

Related to the Anosov-like action conjecture, Baik-Wu-Zhao \cite{baik2024reconstruction} recently defined a flowable group action on the bifoliated plane, and showed that if a group admits a flowable orientation-preserving action on a bifoliated plane, the group is the fundamental group of a 3-manifold with a pseudo-Anosov flow and the bifoliated plane is the orbit space of the flow. The set of axioms is slightly different, but the most important condition in both settings is the hyperbolicity (topologically expanding one leaf and contracting another leaf when two leaves cross). 

\begin{quest} 
Is every Anosov-like action flowable? If not, what are counterexamples and what is the precise difference between two actions? 
\end{quest}

The known examples of groups with a veering pair of invariant laminations are the fundamental groups of 3-manifolds with either a pseudo-Anosov flow without perfect fits or a veering triangulation. The following is a natural question from this. 

\begin{quest} 
Classify all 3-manifolds whose fundamental group admits an action on $\S^1$ with a veering pair of invariant laminations.
\end{quest}

Calegari once conjectured that laminar groups preserving a very full lamination (a lamination all of whose gaps are finite-sided ideal polygons) should be residually finite. The moral of this conjecture is that preserving a lamination with enough structure should restrict the possible algebraic structure of the group. 

\begin{quest}[Calegari] 
Show that a group acts on the circle with a very full invariant lamination is residually finite. 
\end{quest}

\section{Word length on diffeomorphisms groups}
\subsection{Commutator length}
\subsubsection*{Frédéric Le Roux}

Given a group $G$ generated by $\mathcal{S}$, we can defined the \emph{word length}  $\|\cdot \|_{\mathcal{S}}$ by counting the minimal number of generators needed to write a given element, and the \emph{stable word length} of $g$, $\lim \frac{1}{n}\|g^n\|_\mathcal{S}$. When $g$ is perfect and $\mathcal{S}$ is the set of commutators, we get the \emph{commutator length} $cl_\mathcal{S}$ and \emph{stable commutator length} $scl_\mathcal{S}$ on $G$.

A theorem of Michel Herman says that in the group of $\mathrm{Diff}_+^\infty(\S^1)$, every element with diophantine rotation number is $C^\infty$ conjugate to a rotation. Using that the rotation number is a continuous function, it follows easily that every element is the product of a rotation and of a conjugate to a rotation. In particular, every element is the product of two commutators. However, it is not known if one commutator is enough.
\begin{quest}
    Is every element in $\mathrm{Diff}_+^\infty(\S^1)$ a commutator?
\end{quest}

We can ask the same question, say, in $\mathrm{Diff}_0^\infty(\S^2)$.
Here is a variation on $\mathrm{Diff}_0^\infty(\Sigma_g)$, where $\Sigma_g$ is surface of genus $g \geq 1$. After the work of Bowden, Hensel and Webb using their fine curve graph \cite{bowden2022quasi}, we know that the stable commutator length is unbounded on these groups.

\begin{quest}\label{Leroux}
    Can we find $f \in \mathrm{Diff}_0^\infty(\Sigma_g)$ such that $\mathrm{scl}(f) = 0$ but $f$ is not a commutator? What about $\mathrm{Homeo}_0(\Sigma_g)$?
\end{quest}
Somehow, these questions ask for some tools allowing to prove that an element is not a commutator without using quasi-morphisms...

Likewise, on these groups we may define the \emph{flow length} as the word length for the set $\mathcal{S}$ of all time one map of flows, which generates $\mathrm{Homeo}_0(M)$ or $\mathrm{Diff}_0^\infty(M)$ since these groups are simple.
This norm (or its asymptotic version) measures in some sense the dynamical complexity of a map, especially on surfaces where flows are considered to have ``simple" dynamics.
The methods used in \cite{bowden2022quasi} for finding elements with positive scl also yield positive stable flow length. Also note that an element which is conjugate to its inverse has vanishing scl.
\begin{quest}
    Can we find an element in $\mathrm{Homeo}_0(\Sigma_g)$ with vanishing scl and positive stable flow length? Does there exist such an element which is conjugate to its inverse?
\end{quest}
\subsubsection*{Sam Nariman} Question \ref{Leroux} has been answered negatively by Theorem 4.1 in Matsumoto's paper \cite{matsumoto1987some}.

\subsection{Commutators of area preserving homeomorphisms of $\S^2$}
\subsubsection*{Vincent Humili\`ere}

Let $G$ denote the group of area and orientation preserving homeomorphisms of the 2-sphere, and $[G,G]$ its subgroup generated by commutators. We recently proved (\cite{cristofaro2023pfh}, see also \cite{cristofaro2024proof, cristifaro-gardiner2022quantitative}) that $G$ is not perfect, i.e.
\[ [G,G]\varsubsetneq G\]
However, the group $G$ remains poorly understood and the following questions are wide open:

\begin{quest}
    What is $[G,G]$ ? 
 What is the quotient $G/[G,G]$?
\end{quest}

The problem of the first question is to find an explicit description of those elements which are products of commutators. For example, in the group of area-preserving smooth diffeomorphisms of a closed surface $\Sigma$, it is known \cite{banyaga1978structure} that the commutator subgroup consists exactly in those elements in the connected components of identity which have vanishing rotation vector in $H_1(\Sigma;\R)$. It would be great to have a similar "easy" description for our group $G$. However, unfortunately, we cannot hope for a too simple description since it is known that $[G, G]$ is dense in $G$.

For the second question on the quotient $Q=G/[G, G]$, the only thing we know is that there exists an injective group homomorphism $\R\to Q$. We do not know whether $Q$ has torsion elements, nor whether it is a divisible group.

\subsection{Distortion elements and conjugacy classes of elements}
\subsubsection*{Emmanuel Militon}
In a finitely generated group $G$ with finite generating set $S$, we can define the word length $l_S(g)$ of an element $g$ of $G$ as the minimal number of factors you need to write $g$ as a product of an element of $S \cup S^{-1}$.

In such a group, an element $g \in G$ is said to be \emph{distorted} if $$\lim_{n \rightarrow +\infty} \frac{l_S(g^n)}{n}=0.$$
The above limit always exists and this definition is independent of the chosen generating set $S$.

For instance, the element which generates the center of the Heisenberg group with integral coefficients is distorted. As a consequence, any higher rank nonuniform lattice (like $\mathrm{SL}_n(\mathbb{Z})$, for $n \geq 3$) contains plenty of distorted elements.

In any group, we say that an element is distorted if it is contained in a finitely generated subgroup in which it is distorted.

Using such distorted elements, Franks and Handel proved in \cite{franks2006distorsion} that any group morphism from such a higher rank lattice to the group of $C^1$-diffeomorphisms of a higher genus surface which preserves an area form is almost trivial, in the sense that it has finite image.

Since then, distortion elements in groups of homeomorphisms and diffeomorphisms of manifolds have been studied in many articles. For instance, Calegari and Freedman proved in \cite{calegari2006distortion} that all elements of the group of homeomorphisms of a sphere (in any dimension) are distorted and Avila proved that any irrational rotation is distorted in the group of $C^{\infty}$-diffeomorphisms of the circle (see \cite{avila2008distortion}). 

Here we are interested in the group $\mathrm{Homeo}_0(\Sigma)$ which is the identity component of the group of homeomorphisms of an oriented surface $\Sigma$. We assume that the genus of $\Sigma$ is at least $1$. The following obstruction to being distorted was known to Franks and Handel.

Let $\tilde{\Sigma}$ be the universal cover of $\Sigma$ and fix a fundamental domain $D \subset \tilde{\Sigma}$ of the covering map $\tilde{\Sigma} \rightarrow \Sigma$. Also, fix a distance $d$ on $\Sigma$ (compatible with the topology) and lift it as a distance on $\tilde{\Sigma}$ which is invariant under the group of deck transformations.

\begin{prop}
If $f$ is distorted in $\mathrm{Homeo}_0(\Sigma)$, then
$$\lim_{n \rightarrow +\infty} \frac{\mathrm{diameter}(\tilde{f}^{n}(D))}{n}=0.$$ 
\end{prop}
We will call an element that satisfies the conclusion of the proposition non-spreading.

\begin{quest} Is the converse true?
\end{quest}

Setting $d_n := \mathrm{diameter}(\tilde{f}^{n}(D))$, it is proved in \cite{militon2014distortion} that if
$$\lim_{n \rightarrow +\infty} \frac{d_n \log d_n}{n}=0,$$
then $f$ is indeed distorted. 

In the article \cite{militon2018conjugacy}, it is proved that if $f$ has conjugate arbitrarily close to a distorted elements then this element is itself a distorted element. Hence, the above question reduces to the following question (as it is known that translation of the $2$-torus are distorted).

\begin{quest}
Is it true that a non-spreading homeomorphism of the $2$-torus has conjugates arbitrarily close to a rotation ? \\
Is it true that a non-spreading homeomorphism of a higher genus surface has conjugates arbitrarily close to the identity?
\end{quest}

\subsection{Distortion elements in the group of germs}

\subsubsection*{H. Eynard-Bontemps \& A. Navas.}
Let's consider the group $\mathcal{G}^r$ of germs of $C^r$ diffeomorphisms about the origin in $\R$, where $r$ is at least $1$ (the reason for excluding $r=0$ is that, in the continuous setting, every 1-dimensional germ -or homeomorphism- is a distortion element, because every nontrivial element is conjugate to its square, for example). Here, $r$ can be also equal to $\omega$, which means that we also consider the real-analytic case. 

\begin{quest} What are the distortion elements of $\mathcal{G}^r$?
\end{quest}

An example of a distortion element in $\mathcal{G}^{\omega}$ is the germ $$g(x) = \frac{x}{1+x} = x-x^2+x^3-x^4+\ldots.$$ Indeed, if we let $f(x) = x/2$, then we have the relation $fgf^{-1} = g^2$ (actually, $f$ and $g$ generate a copy of the Baumslag-Solitar group $B(1,2)$ inside the affine subgroup of the group of germs). 

A challenging particular case for the question above is that of the germ $g(x) = x-x^2$. Although this germ is $C^1$ conjugate to its square (and, in general, to any nontrivial positive power), it is not $C^2$ conjugate to it, because of the nonvanishing of the Schwarzian derivative at the origin (this is straightforward in class $C^3$, but highly non obvious in class $C^2$ because of the failure of the definition of the Schwarzian derivative in this setting; see \cite{eynard2024residues} on all of this). Thus, in the positive case, proving distortion in this setting would require a more involved idea than just conjugate the germ to a power.

The specific case $r=\omega$ seems to be particular. Indeed, in this framework, there is another conjugacy invariant, the so-called \'Ecalle-Voronin invariant, which is a kind of analog of the Mather invariant. So, it may happen that it should be taken into account in the general problem of distortion in $\mathcal{G}^{\omega}$. The intuition comes from the fact that the Mather invariant plays a major role in the problem of distortion for $C^2$ parabolic diffeomorphisms of the interval because of its direct relation with the asymptotic variation, whose positivity is an obvious obstruction to distortion.

\section{Real analytic foliations}

\subsection{Codimension 1 real analytic foliations}
\subsubsection*{Victor Kleptsyn}
For a finitely generated group of analytic diffeomorphisms of the circle that acts with a Cantor minimal set, a theorem by \'E.~Ghys~\cite{ghys1987classe} states that this group is virtually free. This was used in a work of B.~Deroin, A.~Navas and myself~\cite{deroin2018ergodic} to establish the dynamical properties of such action, in particular, establishing that the Cantor set is necessarily of Lebesgue measure zero. This work was later generalised~\cite{alvarez2019groups} to the case of analytic actions of groups with infinitely many ends, then leading in~\cite{alonso2023ping,alvarez2023ping} to the proof of the Dippolito conjecture on the semi-conjugacy to an action by piecewise affine maps.

For the case of real-analytic codimension 1 foliations, the presence of an exceptional minimal set implies that an extremal leaf (passing through an endpoint of an interval of its complement in a transversal section) has infinitely many ends. 

\begin{quest}
    Can the above theory for the analytic group actions be generalized to the case of codimension~1 real analytic foliations? Can it lead to establishing of the Dippolito conjecture for foliations? 
\end{quest}

\subsection{Exceptional minimal sets}
 \subsubsection*{Steve Hurder}
Let $\mathcal F$ be a $C^r$-foliation of a closed manifold $M$ of codimension-$q$, for $r \geq 1$. Let ${\mathfrak M} \subset M$ be a minimal set for $\mathcal F$. That is, ${\mathfrak M}$ is a closed subset which is a union of leaves, and every leaf in ${\mathfrak M}$ is dense. We say that ${\mathfrak M}$ is \emph{exceptional} if its intersection with every transversal to $\mathcal F$ is totally disconnected. 

The construction of exceptional minimal sets for codimension-$1$ foliations by Sacksteder \cite{sacksteder1964existence} and Rosenberg and Roussarie \cite{rosenberg1970feuilles}, along with Sacksteder's theorem that the Denjoy minimal set does not exist for $C^2$-foliations \cite{sacksteder1965foliations}, can be viewed as the beginnings of the modern study of foliation dynamics.  The works \cite{cantwell1988foliations,cantwell2002endsets,hurder1991exceptional,inaba1988examples,inaba1990resilient,matsumoto1988measure} all consider the structure of exceptional minimal sets in codimension-$1$. The more recent works of Deroin, Kleptsyn and Navas \cite{deroin2009question,deroin2017paradigm} have made great strides towards understanding the dynamics of minimal sets in codimension 1, especially for the case of real analytic foliations.

The analogous questions for higher codimension foliations are essentially completely open.  
\begin{quest}
Classify the exceptional minimal sets for codimension-$q$ foliations of compact manifolds, for $q \geq 2$ and $r \geq 2$.
\end{quest} 

 In codimension-$1$, there is a dichotomy for the derivative of the holonomy maps on a minimal set, in that they either have slow (subexponential) growth, or exponential growth \cite{sacksteder1964existence,hurder1991exceptional,hurder2011lectures}. For higher codimensions, the transverse dynamics can be more complicated, as evidenced by the constructions of solenoidal minimal sets for flows of 3-manifolds, and hyperbolic minimal sets which have the transverse structure of a horseshoe. The work of Clark and Hurder \cite{clark2011embedding} constructed solenoidal minimal sets for foliations with leaf dimensions at least 2, and this work has many references to prior works. A basic problem remains finding new constructions for such sets.
 
 On the other hand, the work on the structure of laminations arising from tiling spaces has led to new topological classification results for exceptional minimal sets, and part of the question is to what extent does the fact that ${\mathfrak M}$ arises from a $C^2$-dynamical system restrict its topological type?

\subsection{Powers of Euler class of flat circle bundles}
\subsubsection*{Yoshihiko Mitsumatsu}

Morita (1984) showed by using the Mather-Thurston map for 
flat $\S^1$-products and the rational homotopy theory 
of $\Lambda B\overline{\Gamma}^\infty_1$ 
that $\chi_{\mathbb{Q}}^k\ne 0 \in 
H^{2k}(B\text{Diff}_+^\infty(\S^1)^{\delta}; {\mathbb{Q}})$ 
for any $k\in{\mathbb{N}}$, where 
$\chi_{\mathbb{Q}}$ is the universal Euler class of flat $\S^1$-bundles 
with rational coefficients. 
Soon later, Haefliger pointed out 
that without the rational homotopy theory but using 
the Mather-Thurston theorem, it can be proved rather easily 
and it is established by Nariman \cite{MR3737509}. 
Also Ghys and Sergiescu \cite{MR896095} gave a different proof, 
and Kitano-Mitsumatsu-Morita gave also yet another
proof (but close to Haefliger-Nariman's one), while 
every proof uses Mather-Thurston theory.   
The nontrivial example for the Euler class $\chi$ itself 
is easily given by closed Riemann surfaces of genus greater than one,  
while no examples are known for the powers. 

\begin{quest}
Give an example of $\chi^2\ne 0$, namely, 
construct a flat $\S^1$-bundle over some oriented closed 4-manifold $M$ 
with $\chi^2\ne 0 \in H^4(M;{\mathbb{Z}})$.  
\end{quest}

Nothing is known about the power of the Euler class 
for real analytic flat $\S^1$-bundles. 

\begin{quest}
$\chi_{\mathbb{Q}}^k\ne 0 \in 
H^{2k}(B\textnormal{Diff}_+^\omega(\S^1)^{\delta}; {\mathbb{Q}})$ for $k\geq 2$ ? 
\end{quest}

There seem to exist very few reasons for that 
the Mather-Thurston theorem holds for `real analytic' 
flat bundles.  Recently we showed (with Morita and Kitano) that 
for flat $\S^1$-bundles not only homologically 
but on the level of classifying maps, 
the Mather-Thurston map admits a homotopy left-inverse.  

\subsection{A subgroup $G$ of $\mathrm{Diff}^\omega (\S^1)_0$}
\subsubsection*{Takashi Tsuboi}
Let $G_1$ be the group $\text{PSL}(2, \R) = SL(2, \R)/\Z_2$. For a positive integer $n$, let $p_n: G_n\to G_1$
be the $n$-fold cyclic covering; then $G_n$ has a topological group structure such that $p_n$ is a
homomorphism. 

Consider the subgroup $R_1$ of $G_1$ which consists of rotations. For $n\in \N$, we put
$R_n=p_n^{-1}(R_1)$. Then as a topological group, $R_1$ is isomorphic to $\R/\Z$;
\[
R_1\cong R_1\cong\cdots \cong R_ n \cdots\cong \R/\Z
\]
We fix these isomorphisms. For $n\in \N$, let $G^{(n)}$  be the free product of $G_i$ for $(i = 1 ..... n)$ with amalgamation $R_1\cong R_1\cong\cdots \cong R_ n$.  Then we have a direct system of groups;
$G^{(0)}\hookrightarrow G^{(1)} \hookrightarrow G^{(2)}\hookrightarrow \cdots$. We define $G$ to be the direct limit group; $G:=\varinjlim G^{(n)}$

\begin{quest}
Can we show that $G$ is dense in $\mathrm{Diff}^\omega (\S^1)_0$ ?
\end{quest}
One can show that if 2-cycle of $G$ has 0 Godbillon-Vey invariant, then it is smoothly a boundary of a codimension 1 foliation of a compact 4 manifold.
\\
Just after posing this question, Michele Triestino told me that there is a paper by James Giblin and Vladimir Markovic \cite{giblin2006classification}, which solved this question.

\subsubsection*{Michele Triestino}
The answer to this question should be positive. For the $C^0$ topology, one can directly apply a result of Giblin-Markovic \cite{giblin2006classification}: closed transitive subgroups of $\mathrm{Homeo}(\S^1)$ containing a non-constant continuous path are conjugate to one of the following $\mathrm{SO}(2)$, $\mathrm{PSL}(2)^{(k)}$, $\mathrm{Homeo}(\S^1)^{(k)}$. For $C^\omega$ topology, a possible way to conclude is to look at the algebra generated by the vector fields tangent to the $\mathrm{PSL}(2)^{(k)}$. This is dense in $\mathrm{Vect}(\S^1)$ because of the classification of Lie algebras of 1d vector fields (is this a result of Lie/Cartan?). This should more generally prove that any amalgamated product $\mathrm{PSL}(2)^{(k)}*_{\mathrm{SO}(2)}\mathrm{PSL}(2)^{(l)}$ is dense in some appropriate $\mathrm{Diff}^{\omega}(\S^1)^{(n)}$. Otherwise, perhaps the Cauchy estimates allow to prove $C^\omega$ directly from $C^0$ density?

In relation to this question, I would like to recall another problem stated by Tsuboi:

\begin{quest}
    Is the subgroup $\langle \mathrm{PSL}(2)^{(k)},\mathrm{PSL}(2)^{(l)}\rangle$ in $\mathrm{Diff}^{\omega}(\S^1)$ isomorphic to the amalgamated product $\mathrm{PSL}(2)^{(k)}*_{\mathrm{SO}(2)}\mathrm{PSL}(2)^{(l)}$?
\end{quest}

\section{Homotopy theory of foliations}
\subsection{Some open questions about the classifying space $B\Gamma$.}
\subsubsection*{Ga\"el Meigniez}

The \emph{homotopy theory of foliations} aims at Reeb's ``fundamental question"
to which A. Haefliger
gave its modern form. To fix ideas, let us restrict
ourselves to the smooth ($C^\infty$) case.
 \emph{Given a $p$-plane field
 on a $n$-manifold, 
under which conditions is it homotopic to a smooth integrable one?}

Some impressive results have been obtained in the golden Lustrum
1969-1976, in particular by W. Thurston; and a few consistent ones since; but
some major questions remain open. Today, the revival
of the theory has a geometric-topological face using the philosophy and methods
 of Gromov's ``h-principle"; and a homotopy-theoretic face.
 (See Nariman \cite{nariman2024foliations} for a recent review).

A central notion is \emph{Haefliger structures,}
also known as $\Gamma$-structures. Vaguely, a $\Gamma_q$-structure
is a kind of singular foliation of codimension $q$.
   Precisely and
   concretely, on a manifold $M$, such a structure amounts to a \emph{microfoliation:} a real vector bundle $\pi:E\to M$ of rank $q$
   over $M$ (the ``normal bundle" of the structure)
    together with a germ, along the zero section
   in $E$,
   of foliation of codimension $q$
   transverse to the fibres of $\pi$. A $\bar\Gamma_q$-structure
   is a $\Gamma_q$-structure whose normal bundle is $M\times\R^q$.

 The central object
  is, for every integer
   $q\ge 0$, Haefliger's space
    $B\Gamma_q$ classifying  the $\Gamma_q$-structures.
    It fibres over $BO(q)$ with fibre $B\bar\Gamma_q$.
       (See
    Tsuboi \cite{tsuboi2009classifying} for an introductory survey on $B\Gamma$).

A major feature of the theory is
a close relation between $B\Gamma$
 and the homology of the
diffeomorphism groups (made discrete),
 through the ``Mather-Thurston theorem".
 
The following questions, which
 seem to resist to the contemporary attempts,
all ask for some triviality in the homotopy of
$B\Gamma$, that would imply the existence of foliations,
as well as some acyclicity in the homology
of some diffeomorphism groups.
On the other hand, if the
answers were negative,
 it would mean that Reeb's fundamental problem admits obstructions
of a completely new nature, different
from the known ones: essentially,
 Bott's obstructions, and
``secondary characteristic classes" such as the Godbillon-Vey
invariant.

\begin{quest}[``Haefliger-Thurston conjecture"]
    \label{A_quest}\
    
 Is $B\bar\Gamma_q$ $2q$-con\-nected i.e. do $\pi_i(B\bar\Gamma_q)$ vanish for all $i\leq 2q$?
 \end{quest}

   Thurston proved that $B\bar\Gamma_q$ is ($q+1$)-connected, and not
   ($2q+1$)-connected
   \cite{thurston1974foliations}. Haefliger and Bott proved that the continuous
    cohomology of $B\bar\Gamma_q$, in a certain sense, vanishes up to degree $2q$ \cite{haefliger1979differential}.  Question \ref{A_quest},
    which arose in 1972 at IAS Princeton during the collective discussions of
   the special year on Foliations \cite{lawson},
    was promoted to a conjecture by Thurston \cite{thurston1974foliations};
     regarded as a question by Lawson \cite{lawson1974foliations}; while Haefliger wrote
   ``one is tempted to conjecture that..." \cite{haefliger1979differential}.
   
   The question if $\pi_{p+q}(B\bar\Gamma_q)$
   vanishes remains open for every pair $(p,q)$ such that
    $2\le p\le q$.  Does every
   $\bar\Gamma_2$-structure on $\S^4$ extend to $\D^5$?
   
As for Haefliger structures
in nonsmooth regularity classes,
Tsuboi obtained results analogous to Question \ref{A_quest}
 in low differentiability classes
\cite{tsuboi1985homology,tsuboi1989foliated}. 
In the piecewise linear
   class, recently Nariman,
    using some model for $B\bar\Gamma_2^{PL}$ due to Greenberg, proved that $B\bar\Gamma_2^{PL}$
   is $4$-connected \cite{nariman2024pl}.
   
   A positive answer to the Haefliger-Thurston conjecture
    would have large-scope
   consequences. The answer to
   Reebs's fundamental question above would always be positive provided that $p\le(n+1)/2$, inducing
   the existence of many foliations~!
    The low-degree homology of several
   diffeomorphism groups (with compact support)
    made discrete would also follow; e.g.
   
 \begin{itemize}
 \item  $H_p(\diff_c(\R^q)^\delta)$
  would vanish  for $p\le q\le 3$;
   
   \item  For  the orientable closed connected surface $\Sigma_g$ 
    of genus $g\ge 2$, 
    the group of the
     orientation-preserving diffeomorphisms
   made discrete would have
     the same degree-$2$ homology as the Mapping class group:
    $$H_2(B{\diff}_+(\Sigma_g)^\delta)
    \cong H_2(B(\pi_0(\diff_+(\Sigma_g))))$$

   \item As for the spheres of dimensions $2$ and $3$,
    $$H_p(B{\diff}_+(\S^q)^\delta)\cong
   H_p(B\textnormal{SO}(q+1))$$ would hold for  $p\le q\le 3$;

 \item  As for hyperbolic $3$-manifolds, one would have
 $$H_p(B{\diff}_+(M)^\delta)
    \cong H_p(B(\textnormal{Isom}_+(M))$$ for $p=2$ and $p=3$,
    for every closed orientable hyperbolic
    $3$-manifold $M$.
    \end{itemize}
 \begin{quest}[Hurder-Vogt \cite{hurder1993product,hurder2005problem,langevin1992list}]\label{B_quest}
 
  Consider the $\bar\Gamma_q$-structure $\gamma$ on the manifold $SO(q)$
 whose microfoliation is the suspension of the canonical
 action of $SO(q)$ on $\R^q$. Is $\gamma$ 
 null-homotopic in $B\bar\Gamma_q$~?
 \end{quest}
  In other words, very concretely, is there a 
  $\bar\Gamma_q$-structure $\bar\gamma$ on $SO(q)\times[0,1]$
  whose restriction over $SO(q)\times 0$ is trivial
  (the microfoliation being the slice foliation of $SO(q)\times\R^q$
   parallel to $SO(q)$) while the restriction
   of $\bar\gamma$ over $SO(q)\times 1$ is
   $\gamma$~?

 The answer is positive for $q\le 4$ (Hurder).
 It is somehow frustrating that we cannot answer for such
 a simple Haefliger structure, showing how few we know!
  If the answer to Question \ref{B_quest} were positive
 for every $q$, then
 the following would hold.
 
 \begin{enumerate}\item $\Omega B\Gamma_q^+$ would
 split homotopically as $$\Omega B\Gamma_q^+=
 \Omega B\bar\Gamma_q\times SO(q)$$
 
 \item The answer to Reeb's fundamental question
 above would always be positive if $M$ is
 homeomorphic to the $n$-sphere;
 
 \item The real cohomology of $BSO(q)$ would inject into
 the real cohomology of $B\Gamma^+_q$, up to degree $2q$.
 \end{enumerate}
 Here $B\Gamma^+_q$ is the classifying space for the
  Haefliger
 structures of
  codimension $q$ with arbitrary, oriented normal bundle.
  Hurder conjectured that (2)
 holds for every $q$; Vogt gave this conjecture
  the very concrete form of Question \ref{B_quest}.
 The would-be corollary (3), by itself a classical question,
  would be a kind of
 reciproque to Bott's obstructions, which state that
 the real cohomology algebra of $BSO(q)$ generated by the Pontrjagin
 classes vanishes in $B\Gamma^+_q$ in all degrees $>2q$.

A probably more accessible question asks if every homotopy
group
 $\pi_p(SO(q)))$
 ($p\ge 0$) vanishes in $\pi_p(B\bar\Gamma_q)$ under the classifying map of $\gamma$.
 
 The Hurder-Vogt question is a particular case of a more general
 and an older question.
 
\begin{quest}[``Foliations with only one leaf": see Tsuboi \cite{tsuboi1984gamma1}]\label{C_quest}\
  Let $\Gamma_q$ be the groupoid of the germs of
  local
 diffeomorphisms of $\R^q$, and $G_q\subset\Gamma_q$
 be the subgroup fixing $0$.  Let $B\bar G_q$ be the
 homotopy-theoretic fibre of the map
 $BG_q\to BGL_q(\R)$ classifying the normal bundle
 of the universal $G_q$-structure.
  Is $B\bar G_q$
 homotopically trivial in $B\bar\Gamma_q$~?
 \end{quest}

 
 This amounts to asking if
  any (framed)
  Haefliger structure whose microfoliation has
 the zero section as a leaf is homotopically trivial.
 Tsuboi answered positively for $q=1$.
 
 \begin{quest}[E. Ghys, appendix to
 \cite{molino1988riemannian}]\label{D_quest}\
 
 i) Is every codimension-$q$ Riemannian
 foliation homotopically trivial in $B\Gamma_q$~?
 
 ii) Is every codimension-$q$ Lie
 foliation homotopically trivial in $B\bar\Gamma_q$~?
 \end{quest}
 Recall that a foliation is \emph{Riemannian}
 is the ambient manifold admits a Riemannian metric for which
 the distance between any two neighbouring leaves is locally
 constant. A special case is, for every Lie group
 $G$, the \emph{$G$-Lie foliations,} whose transverse dynamics
 are modelled on $G$ modulo its left translations.
The classifying space $B\bar G$ for the $G$-Lie foliations
is a principal $G$-bundle over the Eilenberg-McLane space
$K(G^\delta,1)$; the above question amounts to know
if $B\bar G$ is null-homotopic in $B\bar\Gamma_{\dim(G)}$.
As Ghys noticed, using an argument of Tsuboi's,
 the answer is positive if $G$ is $\R^q$;
and probably also if $G$ is nilpotent.

\subsection{Piecewise Linear homeomorphisms}
\subsubsection*{Sam Nariman}
The subgroups of piecewise linear (PL for short) homeomorphisms of the line have been a rich source for interesting finitely generated groups with surprising algebraic and dynamical properties. However, not much is known about the algebraic properties of the PL homeomorphisms of higher dimensional PL manifolds. Although they are more combinatorial in nature, the analytical tools for diffeomorphism groups and the Mather infinite repetition trick for homeomorphisms are not available for PL homeomorphisms. So the following basic question due to Epstein \cite{epstein1970simplicity} is still open.
\begin{quest}(Epstein) Let $M$ be a PL manifold. Is $\textnormal{PL}_0(M)$, the group of PL homeomorphisms of $M$ that are isotopic to the identity, a simple group? 
\end{quest}
By Epstein's result, perfectness and simplicity are equivalent in this case and he proved that $\text{PL}_c(\R)$ and $\text{PL}_0(\S^1)$ are perfect by observing that in dimension one,  PL homeomorphisms are generated by certain ``typical elements'' and those typical elements can be easily written as commutators. One approach to solve Epstein's problem is to prove that every element in $\text{PL}_0(M)$ is generated by ``typical" elements (see \cite{epstein1970simplicity, nariman2024foliations} for definitions).

I took a different approach in dimension $2$ and used Greenberg's work \cite{greenberg1987classifying, greenberg1992generators} on  the Haefliger classifying space of PL foliations, to first show that 
\begin{thm}(\cite{nariman2024pl})
The space $\overline{\mathrm{B}\Gamma}_2^{\textnormal{\text{PL}}}$ is  $4$-connected.
\end{thm}
where $\overline{\mathrm{B}\Gamma}_2^{\textnormal{\text{PL}}}$ is classifying space of codimension $2$ PL Haefliger structures. Using Mather-Thurston's theorem for PL homeomorphisms, one gets as a corollary.
\begin{cor}(\cite[Corollary 1.4]{nariman2024pl})
Let $\Sigma$ be an oriented compact surface, possibly with a boundary. The identity component ${\textnormal{\text{PL}}}_\circ(\Sigma,\text{rel }\partial)$ is a simple group as a discrete group.
\end{cor}
It is likely that using Bowden-Hensel-Webb's work, one could show that ${\textnormal{\text{PL}}}_\circ(\Sigma,\text{rel }\partial)$ is not uniformly perfect  for cases where $\text{Homeo}_0(\Sigma,\text{rel }\partial)$ is not uniformly perfect. This gives an h-principle approach for perfectness of ${\textnormal{\text{PL}}}_\circ(\Sigma,\text{rel }\partial)$. I still do not know an algebraic proof and more importantly do not know how to generalize my approach to higher dimensions. In particular, I do not know the following.
\begin{quest}
  Is  the space $\overline{\mathrm{B}\Gamma}_n^{\textnormal{\text{PL}}}$ is  $n+1$-connected?
\end{quest}
This perspective also suggests problems about homological properties of $\text{PL}_0(M)$ by the way of analogy. For example, the group homology of $\text{PL}_c(\R)$ has been studied by Greenberg, Tsuboi, Ghys and Sergiescu. In particular, there is a discrete version of Godbillon-Vey class for codimension one PL foliations. Ghys and Sergiescu found the following 2-cocycle for $f_1, f_2 \in \text{PL}_0(\S^1)$, let 
\[
\tiny\overline{gv}(f_1, f_2)= \frac12 \sum_{x\in \S^1} \begin{vmatrix}\log{f_2'}(x+0) & \log{(f_1\circ f_2)'}(x+0)\\ \log{f_2'}(x+0)-\log{f_2'}(x-0) &\log{(f_1\circ f_2)'}(x+0)-\log{(f_1\circ f_2)'}(x-0)\end{vmatrix}
\]
where the summation is over finitely many singularities of $f_1$ and $f_2$.
\begin{quest}
Does an analog of discrete $\textnormal{GV}$ exist in higher dimensions? In particular, does there exist a $3$-cocycle in the third group cohomology $H^3(\textnormal{PL}_c(\R^2); \R)$ that induces a surjection
\[
H_3(\textnormal{PL}_c(\R^2); \mathbb{Z})\twoheadrightarrow \R.
\]
\end{quest}

Also, it is tempting to think about the analog of Milnor's conjecture for PL homeomorphisms. Recall that for a finite dimensional Lie group $G$, Milnor conjectured that 
\[
H_*(BG^{\delta}; \mathbb{F}_p)\to H_*(BG; \mathbb{F}_p)
\]
is an isomorphism in all degrees. Here, $\mathbb{F}_p$ is the finite field with $p$ elements, $BG$ is the classifying space of the Lie group $G$ and $BG^{\delta}$ is the classifying space of $G$ made discrete. The work of Greenberg on PL foliations seems to suggest that the classifying space of PL Haefliger structures is inductively built out of classifying spaces of Lie groups made discrete. So it might be reasonable to expect the following.

\begin{quest}
Is the map 
\[
H_*(B\textnormal{PL}(M)^{\delta}; \mathbb{F}_p)\to H_*(B\textnormal{PL}(M); \mathbb{F}_p)
\]
an isomorphism in all degrees?
\end{quest}
The work of Greenberg implies that this is true for $\text{PL}_c(\R)$.

\subsection{Haefliger's classifying spaces $B\overline{\Gamma}^{\infty / \omega}_1$}
\subsubsection{Yoshihiko Mitsumatsu}

For $r=\infty$, of course the $2q$-connectivity (for $q\geq 2$) is
one of the greatest and most important problems.

For $r=\omega$ $q\geq 2$ not very much is known.
For $q=1$, Haefliger proved that $B\overline{\Gamma}^\omega_1$ is
a $K(\pi, 1)$ space for a big discrete group $\pi=\Gamma_H$ which
we call the {\it Haefliger group}. He also showed that
$\Gamma_H$ is perfect, continuously generated, and torsion-free.
Therefore we have
$H_1(B\overline{\Gamma}^\omega_1;{\mathbb{Z}})=H_1(\Gamma_H;{\mathbb{Z}})=0$.  
As the Thurston (or Bott-Brooks) cycles show that 
the Godbillon-Vey invariant gives a surjection 
$H_3(B\Gamma_H; {\mathbb{Z}}) \to {\R}$ 
which factors through 
$H_3(B\overline{\Gamma}^\infty_1; {\mathbb{Z}})
=\pi_3(B\overline{\Gamma}^\infty_1)$.  
Tsuboi \cite{MR744853} constructed a homomorphic section 
${\R} \to \pi_3(B\overline{\Gamma}^\infty_1)$ 
and recently Morita (et al.) modified it into 
${\R} \to H_3(B\Gamma_H;{\mathbb{Z}})$.   
So we have 
$\pi_3(B\overline{\Gamma}^\infty_1)\cong {\R}\oplus \mathrm{ker}GV$ 
as well as 
$H_3(B\Gamma_H;{\mathbb{Z}})
\cong {\R}\oplus \mathrm{ker} GV$. 

\begin{quest}Are these two $\mathrm{ker }GV$'s trivial or not?  
\end{quest}

As $B\Gamma_H$ has a trivial $\pi_2$ there is no nontrivial spherical cycle in
$H_2(B\Gamma_H; {\mathbb{Z}})$.  
It is also known that all maps from the  
2-torus $T^2$ to $B\Gamma_H$  yield trivial cycles in $H_2(B\Gamma_H; {\mathbb{Z}})$.

\begin{quest}(Tsuboi \cite{MR3726712}) Determine whether $H_2(B\Gamma_H;{\mathbb{Z}})=0$.  
\end{quest}

If a loop in $B\Gamma_H$ is freely homotopic to one that is 
transverse to the $\Gamma^\omega_1$ structure, it presents 
a non trivial element in $\pi_1(B\overline{\Gamma}^\omega_1)=\Gamma_H$ 
which is of infinite order. 

\begin{quest}
Does there exist a homotopically non-trivial loop in 
$B\overline{\Gamma}^\omega_1$ which is not free-homotopic to 
a transverse one? 
\end{quest} 

If this is not true, we can describe the whole $\Gamma_H$ fairly clearly, 
but it is not plausible. 
If it is true, we have a ${\mathbb{Z}}_{\geq 0}$-valued pseudonorm 
on $\Gamma_H$.

\begin{quest}
Is the group $\Gamma_H$ simple? 
\end{quest}%

	\bibliographystyle{alpha}
	\bibliography{Reference-CIRM-problem-list.bib}

\end{document}